\documentclass{llncs}
\usepackage[]{graphicx}
\usepackage{color}
\usepackage{amsmath}
\usepackage{mathtools}
\usepackage{amssymb}
\usepackage{pbox}
\usepackage{booktabs}
\usepackage{colonequals}
\usepackage{ulem}

\usepackage{algorithm}
\usepackage[noend]{algpseudocode}
\begin{document}

\title{Entropic Trace Estimates for Log Determinants
}

\titlerunning{Entropic Trace Estimates for Log Determinants}        

\author{Jack Fitzsimons\inst{1}         \and
        Diego Granziol\inst{1}       \and
        Kurt Cutajar\inst{2}         \\
        Michael Osborne\inst{1}        \and
        Maurizio Filippone\inst{2}          \and
        Stephen Roberts\inst{1}
}

\institute{Department of Engineering, University of Oxford, UK\\
\email{\{jack,diego,mosb,sjrob\}@robots.ox.ac.uk}
\and Department of Data Science, EURECOM, France \\
\email{\{kurt.cutajar,maurizio.filippone\}@eurecom.fr}
}



\maketitle

\begin{abstract}
The scalable calculation of matrix determinants has been a bottleneck to the widespread application of many machine learning methods such as determinantal point processes, Gaussian processes, generalised Markov random fields, graph models and many others.
In this work, we estimate log determinants under the framework of maximum entropy, given information in the form of moment constraints from stochastic trace estimation.
The estimates demonstrate a significant improvement on state-of-the-art alternative methods, as shown on a wide variety of UFL sparse matrices. 
By taking the example of a general Markov random field, we also demonstrate how this approach can significantly accelerate inference in large-scale learning methods involving the log determinant.
\end{abstract}

\newcommand{\ce}{\colonequals}
\newcommand{\ec}{\equalscolon}

\definecolor{mygreen}{rgb}{0.2, 0.7, 0.2}
\definecolor{myorange}{rgb}{0.9, 0.5, 0.0}
\definecolor{myamber}{rgb}{1.0, 0.75, 0.0}
\definecolor{mylilac}{rgb}{0.6, 0.4, 0.8}

\newcommand\noteMF[1]{\textcolor{red}{MF - #1}}
\newcommand\noteKC[1]{\textcolor{mylilac}{KC - #1}}
\newcommand\noteJF[1]{\textcolor{myorange}{JF - #1}}
\newcommand\noteMO[1]{\textcolor{mygreen}{MO - #1}}
\newcommand\noteSR[1]{\textcolor{blue}{SR - #1}}
\newcommand\noteDG[1]{\textcolor{myamber}{DG - #1}}

\newcommand\highlight[1]{\textcolor{red}{#1}}

\newcommand{\nobs}{n} 
\newcommand{\R}{\mathbb{R}}
\newcommand{\N}{\mathbb{N}}
\newcommand{\Z}{\mathbb{Z}}
\newcommand{\F}{\mathcal{F}}
\newcommand{\I}{\mathcal{I}}
\newcommand{\LL}{\mathcal{L}}
\newcommand{\uu}{\mathbf{u}}
\newcommand{\ee}{\mathbf{e}}

\newcommand{\E}{\mathrm{E}}
\newcommand{\const}{\mathrm{const.}}
\newcommand{\diag}{\mathrm{diag}}
\newcommand{\Tr}{\mathrm{Tr}}
\newcommand{\Det}{\mathrm{Det}}
\newcommand{\GP}{\mathrm{GP}}
\newcommand{\ud}{\mathrm{d}}

\newcommand{\norm}{\mathcal{N}}

\newcommand{\avect}{\mathbf{a}}
\newcommand{\dvect}{\mathbf{d}}
\newcommand{\fvect}{\mathbf{f}}
\newcommand{\gvect}{\mathbf{g}}
\newcommand{\hvect}{\mathbf{h}}
\newcommand{\mvect}{\mathbf{m}}
\newcommand{\pvect}{\mathbf{p}}
\newcommand{\svect}{\mathbf{s}}
\newcommand{\uvect}{\mathbf{u}}
\newcommand{\vvect}{\mathbf{v}}
\newcommand{\zvect}{\mathbf{z}}
\newcommand{\xvect}{\mathbf{x}}
\newcommand{\yvect}{\mathbf{y}}
\newcommand{\wvect}{\mathbf{w}}
\newcommand{\Wvect}{\mathbf{W}}
\newcommand{\tvect}{\mathbf{t}}
\newcommand{\zerovect}{\mathbf{0}}
\newcommand{\onesvect}{\mathbf{1}}

\newcommand{\betavect}{\boldsymbol{\beta}}
\newcommand{\thetavect}{\boldsymbol{\theta}}
\newcommand{\Thetavect}{\mathbf{\Theta}}
\newcommand{\psivect}{\boldsymbol{\psi}}
\newcommand{\Psivect}{\boldsymbol{\Psi}}
\newcommand{\etavect}{\boldsymbol{\eta}}
\newcommand{\rhovect}{\boldsymbol{\rho}}
\newcommand{\tauvect}{\boldsymbol{\tau}}
\newcommand{\nuvect}{\boldsymbol{\nu}}
\newcommand{\muvect}{\boldsymbol{\mu}}
\newcommand{\omegavect}{\boldsymbol{\omega}}
\newcommand{\Omegavect}{\mathbf{\Omega}}
\newcommand{\sigmavect}{\boldsymbol{\sigma}}
\newcommand{\zetavect}{\boldsymbol{\zeta}}
\newcommand{\varepsilonvect}{\boldsymbol{\epsilon}}
\newcommand{\deltavect}{\boldsymbol{\delta}}

\newcommand{\bigO}{\mathcal{O}}

\newcommand{\name}[1]{{\textsc{#1}}\xspace}

\newcommand{\mcmc}{\name{mcmc}}

\newcommand{\gp}{\name{gp}}
\newcommand{\ard}{\name{ard}}

\newcommand{\relu}{{\textsc{r}}e\name{lu}}

\newcommand{\arc}{\name{arc}}
\newcommand{\rbf}{\name{rbf}}

\newcommand{\mvp}{\name{mvp}}
\newcommand{\mnll}{\name{mnll}}
\newcommand{\rmse}{\name{rmse}}
\newcommand{\nelbo}{\name{nelbo}}


\section{Introduction}
\label{intro}

Scalability is a key concern for machine learning in the big data era, whereby inference schemes are expected to yield optimal results within a constrained computational budget.
Underlying these algorithms, linear algebraic operations with high computational complexity pose a significant bottleneck to scalability, and the log determinant of a matrix~\cite{Bai1997}  falls firmly within this category of operations.
The canonical solution involving Cholesky decomposition~\cite{Golub1996} for a general $n\times n$ positive definite matrix, $A$, entails time complexity of $\mathcal{O}(n^{3})$ and storage requirements of $\mathcal{O}(n^{2})$, which is infeasible for large matrices.
Consequently, this term greatly hinders widespread use of the learning models where it appears, which includes determinantal point processes~\cite{Macchi1975}, Gaussian processes~\cite{Rasmussen2006}, and graph problems~\cite{Wainwright2006}.

The application of kernel machines to vector valued input data has gained considerable attention in recent years, enabling fast linear algebra techniques.
Examples include Gaussian Markov random fields~\cite{Rue2005} and Kronecker-based algebra~\cite{Saatci2011}, while similar computational speed-ups may also be obtained for sparse matrices.
Nonetheless, such structure can only be expected in selected applications, thus limiting the widespread use of such techniques.

In light of this computational constraint, several approximate inference schemes have been developed for estimating the log determinant of a matrix more efficiently.
Generalised approximation schemes frequently build upon iterative stochastic trace estimation techniques \cite{Avron2011}.
This includes polynomial approximations such as Taylor and Chebyshev expansions \cite{Aune2014,Han2015}.
Recent developments shown to outperform the aforementioned approximations include estimating the trace using stochastic Lanczos quadrature \cite{Ubaru2016}, and a probabilistic numerics approach based on Gaussian process inference which incorporates bound information~\cite{Fitzsimons2017}.
The latter technique is particularly significant as it introduces the possibility of quantifying the numerical uncertainty inherent to the approximation.


In this paper, we present an alternative probabilistic approximation of log determinants rooted in information theory, which exploits the relationship between stochastic trace estimation and the moments of a matrix's eigenspectrum. These estimates are used as moment constraints on the probability distribution of eigenvalues.
This is achieved by maximising the entropy of the probability density $p(\lambda)$ given our moment constraints.
In our inference scheme, we circumvent the issue inherent to the Gaussian process approach~\cite{Fitzsimons2017}, whereby positive probability mass may occur in the region of negative densities.
In contrast, our proposed entropic approach implicitly encodes the constraint that densities are necessarily positive.
Given equivalent moment information, we achieve competitive results on matrices obtained from the SuiteSparse Matrix Collection~\cite{Davis11} which consistently outperform competing approximations to the log-determinant \cite{Fitzsimons2017,Boutsidis2015}.\newline

\noindent The most significant contributions of this work are listed below.
\begin{enumerate}
\item We develop a novel approximation to the log-determinant of a matrix which relies on the principle of maximum entropy enhanced with moment constraints derived from stochastic trace estimation.
\item We present the theory motivating the use of maximum entropy for solving this problem, along with insights on why we expect particularly significant improvements over competing techniques for large matrices.
\item We directly compare the performance of our entropic approach to other state-of-the-art approximations to the log-determinant.
This evaluation covers real sparse matrices obtained from the SuiteSparse Matrix Collection~\cite{Davis11}. 
\item Finally, to showcase how the proposed approach may be applied in a practical scenario, we incorporate our approximation within the computation of the log-likelihood term of a Gaussian Markov random field, where we obtain a significant increase in speed.

\end{enumerate}

\subsection{Related Work}

The methodology presented in this work predominantly draws inspiration from the recently introduced probabilistic numerics approach to estimating the log determinant of a matrix~\cite{Fitzsimons2017}.
In that work, the computation of the log determinant is re-interpreted as a probabilistic estimation problem, whereby results obtained from budgeted computations are used to infer accurate estimates for the log determinant.
In particular, within that proposed framework, the eigenvalues of a matrix $A$ are modelled from noisy observations of $\Tr(A^k)$ obtained from stochastic trace estimation~\cite{Avron2011} using the Taylor approximation method.
By modelling such noisy observations using a Gaussian process~\cite{Rasmussen2006}, Bayesian quadrature~\cite{OHagan91} can then be invoked for making predictions on the infinite series of the Taylor expansion, and in turn estimating the log determinant.
Of particular interest is the uncertainty quantification inherent to this approach, which is a notable step forward in the direction of measuring the complete numerical uncertainty associated with approximating large-scale inference models.
The estimates obtained using this Bayesian set-up may be further improved by considering known upper and lower bounds on the value of the log determinant~\cite{Bai1997}.
In this paper, we provide an alternative to this approach by interpreting the observed moments as being constraints on the probability distribution of eigenvalues underlying the computation of the log determinant.
As we shall explore, our novel entropic formulation makes better calibrated prior assumptions than the previous work, and consequently yields superior performance.

More traditional approaches to approximating the log determinant build upon iterative algorithms, and exploit the fact that the log determinant may be rewritten as the trace of the logarithm of the matrix.
This features in both the Chebyshev expansion approximation~\cite{Han2015}, as well as the widely-used Taylor series approximation upon which the aforementioned probabilistic inference approaches are built.
Recently, an approximation to the log determinant using stochastic Lanczos quadrature~\cite{Ubaru2016} has been shown to outperform the aforementioned polynomial approaches, while also providing probabilistic error bounds.
Finally, given that the logarithm of a matrix often appears multiplied by a vector (for example the log likelihood term of a Gaussian process~\cite{Rasmussen2006}), the spline approximation proposed in~\cite{Chen11} may be used to accelerate computations.

\section{Background}

In this section, we shall formally introduce the concepts underlying the proposed maximum entropy approach to approximating the log determinant.
We start by describing stochastic trace estimation and demonstrate how this can be applied to estimating the trace term of matrix powers.
Subsequently, we illustrate how the latter terms correspond to the raw moments of the matrix's eigenspectrum, and show how the log determinant may be inferred from the distribution of eigenvalues constrained by such moments.

\subsection{Stochastic Trace Estimation}
\label{stoctrace}

Estimating the trace of implicit matrices is a central component of many approaches to approximate the log determinant of a matrix. 
Stochastic trace estimation~\cite{Avron2011} builds a Monte Carlo estimate of the trace of a matrix $A$ by repeatedly multiplying it by \textit{probing vectors} $\zvect$,

\[
\text{Tr}(A) \approx \frac{1}{m} \sum_{i=1}^m \zvect_i^T A \zvect_i,
\]

\noindent such that the expectation of $\zvect_i \zvect_i^T$ is the identity, $\mathbb{E}[\zvect_i \zvect_i^T] = I$. This can be readily seen using the expectation of $\text{Tr}(\zvect_i^T A \zvect_i)$ by exploiting the cyclical property of the trace operation.
As such, many choices of how to sample the probing vectors have emerged.
Possibly the most na\"ive choice involves sampling from the columns of the identity matrix; however, due to poor expected sample variance this is not widely used in the literature. 
Sampling from vectors on the unit hyper-sphere, and correspondingly sampling normal random vectors (Gaussian estimator), significantly reduces the sample variance, but more random bits are required to generate each sample. 
A major progression for stochastic trace estimation was the introduction of Hutchinson's method \cite{Hutchinson1990}, which sampled each element as a Bernoulli random variable requiring only a linear number of random bits, while also reducing the sample variance even further. 
A more recent approach involves sampling from sets of mutually unbiased bases (MUBs)~\cite{Fitzsimons2017}, requiring only a logarithmic number of bits.
Table \ref{tab:comp} (adapted from~\cite{Fitzsimons2017}) provides a concise overview of the landscape of probing vectors.

\renewcommand{\arraystretch}{1.5}
\begin{table*}[tbh]
\caption{Comparison of single shot variance $V$, worst case single shot variance $V^\text{worst}$ and number of random bits $R$ required for commonly used trace estimators and the MUBs estimator. ($^*$ \textit{required for floating point precision}) \label{tab:comp}}
\begin{tabular}{lccc}
\toprule
& $V$ & $V^\text{worst}$ & R\\ 
\midrule
Fixed basis estimator &  $d \sum_{i=1}^d M_{ii}^2 - \text{Tr}(A)^2$   &    $(d-1)\text{Tr}(A)^2$              &$\log_2(d)$\\
MUBs estimator        & $\frac{d}{d+1} \text{Tr}(A^2) - \frac{1}{d+1} \text{Tr}(A)^2$    & $\frac{d-1}{d+1}\text{Tr}(A^2)$ &$2 \log_2(d)$\\
Hutchinson estimator  & $2\left(\text{Tr}(A^2) - \sum_{i=1}^d A_{ii}^2\right)$    & $\frac{2(d-1)}{d} \text{Tr}(A^2)$           & $d$\\
Gaussian estimator & $2\text{Tr}(A^2)$ & $2\text{Tr}(A^2) $& $\mathcal{O}(d)$$^*$  \\
\bottomrule
\end{tabular}
\centering
\end{table*}

A notable application of stochastic trace estimation is the approximation of the trace term for matrix powers, $\text{Tr}(A^k)$. 
Stochastic trace estimation enables vector-matrix multiplications to be propagated right to left, costing $\mathcal{O}(n^2)$, rather than the $\mathcal{O}(n^3)$ complexity required by matrix multiplication.
This simple trick has been applied in several domains such as counting the number of triangles in graphs \cite{avron2010counting}, string pattern matching \cite{atallah2013lower} and of course estimating the log determinant of matrices, as discussed in this work. 

\subsection{Raw Moments of the Eigenspectrum}

The relation between the raw moments of the eigenvalue distribution and the trace of matrix powers allows us to exploit stochastic trace estimation for estimating the log determinant. 
Raw moments are defined as the mean of the random variables raised to integer powers.
Given that the function of a matrix is implicitly applied to its eigenvalues, in the case of matrix powers this corresponds to raising the eigenvalues to a given power.
For example, the $k^{\text{th}}$ raw moment of the  distribution over the eigenvalues (a mixture of Dirac delta functions) is $\sum_{i=1}^m \lambda^k p(\lambda)$, where $p(\lambda)$ is the distribution of eigenvalues.
The first few raw moments of the eigenvalues are trivial to compute. 
Denoting the $k$th raw moment as $\mathbb{E}[\lambda^k]$, we have that $\mathbb{E}[\lambda^0] = 1$, $\mathbb{E}[\lambda^1] = \frac{1}{n} \text{Tr}(A)$ and  $\mathbb{E}[\lambda^2] = \frac{1}{n} \sum_{i,j} A_{i,j}^2$.
More generally, the  $k^{\text{th}}$ raw moment can be formulated as $\mathbb{E}[\lambda^k] = \frac{1}{n} \text{Tr}(A^k)$, which can be estimated using stochastic trace estimation.
These identities can be easily derived using the definitions and well known identities of the trace term and Frobenius norm.

\subsection{Approximating the Log Determinant}

In view of the relation presented in the previous subsection, we can reformulate the log determinant of a matrix in terms of its eigenvalues using the following derivation:
\begin{equation}
 \log\bigl(\Det(A)\bigr) = \sum_{i=1}^{n}\log(\lambda_{i}) \ce n \mathbb{E}\left[\log (\lambda)\right] \approx n\int p(\lambda) \log (\lambda)\ud\lambda
\end{equation}
where the approximation is introduced due to our estimation of $p(\lambda)$, the probability distribution of eigenvalues. If we knew the true distribution of $p(\lambda)$ it would hold with equality.

Given that we can obtain information about the moments of $p(\lambda)$ through stochastic trace estimation, we can solve this integral by employing the principle of maximum entropy, while treating the estimated moments as constraints.
While not explored in this work, it is worth noting that in the event of moment information combined with samples of eigenvalues, we would use the method of maximum relative entropy with data constraints, which is in turn a generalisation of Bayes' rule~\cite{Caticha2012}. This can be applied, for example, in the quantum linear algebraic setting \cite{nielsen2002quantum}. 

\section{Estimating the Log Determinant using Maximum Entropy}
\label{maxent}
The maximum entropy method (MaxEnt) \cite{Presse2013} is a procedure for generating the most conservatively uncertain estimate of a probability distribution possible with the given information, which is particularly valued for being maximally non-committal with regard to missing information \cite{Jaynes1957}.
In particular, to determine a probability density $p(\xvect)$, this corresponds to maximising the functional 
\begin{equation}
\label{BSG}
S = - \int p(\vec{x})\log p(\vec{x})d\vec{x}- \sum_{i}\alpha_{i}\bigg[\int p(\vec{x})f_{i}(\vec{x})d\vec{x} - \mu_{i}\bigg]
\end{equation}
with respect to $p(\vec{x})$, where $ \mathbb{E}[ f_{i}(\vec{x}) ]= \mu_{i}$ are given constraints on the probability density.
The first term in the above equation is referred to as the Boltzmann-Shannon-Gibbs (BSG) entropy, which has been applied in multiple fields, ranging from condensed matter physics~\cite{Giffin2016} to finance~\cite{Neri2012,Buchen1996}.
Along with its path equivalent, maximum caliber \cite{Gonzalez2014}, it has been successfully used to  derive statistical mechanics \cite{Granziol2017}, non-relativistic quantum mechanics, Newton's laws and Bayes' rule \cite{Gonzalez2014,Caticha2012}.
Under the axioms of consistency, uniqueness, invariance under coordinate transformations, sub-set and system independence, it can be proved that for constraints in the form of expected values, drawing self-consistent inferences requires maximising the entropy  \cite{Shore1980,Presse2013}.
Crucial for our investigation are the functional forms $f_{i}(\vec{x})$ of constraints for which the method of maximum entropy is appropriate.
The axioms of Johnson and Shore~\cite{Shore1980} assert that the entropy must have a unique maximum and that the BSG entropy is convex. 
The entropy hence has a unique maximum provided that the constraints are convex.
This is satisfied for any polynomial in $\xvect$ and hence, maximising the entropy given moment constraints constitutes a self-consistent inference scheme  \cite{Presse2013}.

\subsection{Implementation}

Our implementation of the system follows straight from stochastic trace estimation to estimate the raw moments of the eigenvalues, maximum entropy distribution given these moments and, finally, determining the log of the geometric mean of this distribution. 
The log geometric mean is an estimate of the log determinant divided by the dimensionality of $A \in \mathbb{R}^{n \times n}$. 
We explicitly step through the subtleties of the implementation in order to guide the reader through the full procedure.

By taking the partial derivatives of $S$ from Equation \eqref{BSG}, it is possible to show that the maximum entropy distribution given moment information is of the form
\[
p(\lambda) = \exp(-1 + \sum_i \alpha_i \mu_i).
\]
The goal is to find the set of $\alpha_i$ which match the raw moments of $p(\lambda)$ to the observed moments $\{\mu_i\}$. 
While this \textit{may} be performed symbolically, this becomes intractable for larger number of moments, and our experience with current symbolic libraries~\cite{SymPy2017,Mathematica} is that they are not extendable beyond more than 3 moments.
Instead, we turn our attention to numerical optimisation.
Early approaches to optimising maximum entropy coefficients worked well for a small number of coefficients but became highly unstable as the number of observed moments grew \cite{mohammad1992matlab}. 
However, building on these concepts more stable approaches emerged \cite{bandyopadhyay2005maximum}. 
Algorithm \ref{alg:coefopt} outlines a stable approach to this optimisation under the conditions that $\lambda_i$ is strictly positive and the moments lie between zero and one. 
We can satisfy these conditions by normalising our positive definite matrix by the maximum of the Gershgorin intervals~\cite{Gershgorin1931}.

\algrenewcommand\algorithmicindent{1.3em}
\renewcommand{\algorithmicrequire}{\textbf{Input:}}
\renewcommand{\algorithmicensure}{\textbf{Output:}}

\begin{algorithm}
\caption{Optimising the Coefficients of the MaxEnt Distribution}\label{alg:coefopt}
\begin{algorithmic}[1]
\vspace{0.5em}
\Require Moments $\{\mu_i\}$, Tolerance $\epsilon$ 
\Ensure Coefficients $\{\alpha_i\}$
\State $\alpha_i \sim \mathcal{N}(0,1)$
\State $i \gets 0$
\State $p(\lambda) \gets \exp(-1 - \sum_k \alpha_k \lambda^k)$
\While{error $< \epsilon$}
\State $\delta \gets \log \left(\frac{\mu_i}{\int \lambda^i p(\lambda) d\lambda} \right)$
\State $\alpha_i \gets \alpha_i + \delta$
\State $p(\lambda) \gets p(\lambda | \alpha)$
\State error $\gets \max |\int \lambda^i p(\lambda) d\lambda - \mu_i|$
\State $i \gets \text{mod}(i+1, \text{length}(\mu))$
\EndWhile 
\end{algorithmic}
\end{algorithm}

Given Algorithm \ref{alg:coefopt}, the pipeline of our approach can be pieced together.
First, the raw moments of the eigenvalues are estimated using stochastic trace estimation.
These moments are then passed to the maximum entropy optimisation algorithm to produce an estimate of the distribution of eigenvalues, $p(\lambda)$.
Finally, $p(\lambda)$ is used to estimate the log geometric mean of the distribution, $\int \log(\lambda) p(\lambda) d\lambda$.
This term is multiplied by the dimensionality of the matrix and if the matrix is normalised, the log of this normalisation term is added again.
These steps are laid out more concisely in Algorithm \ref{alg:logdet}.  

\algrenewcommand\algorithmicindent{1.3em}
\renewcommand{\algorithmicrequire}{\textbf{Input:}}
\renewcommand{\algorithmicensure}{\textbf{Output:}}

\begin{algorithm}
\caption{Entropic Trace Estimation for Log Determinants}\label{alg:logdet}
\begin{algorithmic}[1]
\vspace{0.5em}
\Require PD Symmetric Matrix $A$, Order of stochastic trace estimation $k$, Tolerance $\epsilon$
\Ensure Log Determinant Approximation $\log|A|$
\State $B = A/\|A\|_2$
\State $\mu$ (moments)$ \gets$ StochasticTraceEstimation$(B, k)$ 
\State $\alpha$ (coefficients) $\gets \text{MaxEntOpt(}\mu, \epsilon)$
\State $p(\lambda) \gets p(\lambda | \alpha)$
\State $\log|A| \gets n\int \log(\lambda) p(\lambda) d\lambda + n\log(\|A\|_2)$
\end{algorithmic}
\end{algorithm}

\section{Insights for Large Matrices}

The method of entropic trace estimation has the interesting property where we expect the relative error to decrease as the matrix size $N$ increases.
Colloquially, we can liken maximum entropy to a maximum likelihood over distributions, where this likelihood functional is raised to the number of particles in the system, corresponding to the number of eigenvalues in the matrix.
Given that there is a global maximum, as the number of eigenvalues increases, the functional tends to a delta functional around the $p(x)$ of maximum entropy.
This confirms that within the scope of our problem's continuous distribution over eigenvalues, whenever the number of eigenvalues (and correspondingly the dimensionality of the matrix) tends towards infinity, we expect the maximum entropy solution to converge to the true solution.
This gives further credence to the suitability of our method when applied to large matrices.
We substantiate this claim by delving into the fundamentals of maximum entropy for physical systems and extending the analogy to functionals over the space of densities.
We prove that in the limit of $N \rightarrow \infty$ that the maximum entropy distribution dominates the space of solutions satisfying the constraints.
We demonstrate the practical significance of this assertion by setting up an experiment using synthetically constructed random matrices, where this is in fact verified.

\subsection{Law of Large Numbers for Maximum Entropy in Matrix Methods}

In order to demonstrate our result, we consider the quantity $W$, which represents the number of ways in which our candidate probability distribution can recreate the observed moment information.
In order to make this quantity finite we consider the discrete distribution characterised by machine precision $\epsilon$.
We show that $W = \exp(NS)$, where S is the entropy.
Hence maximising the entropy is equivalent to maximising $W$, as $N$ is fixed.
In the continuous limit, we consider the ratio of two such terms $F_{i} = W_{i}/\sum_{j}W_{j}$, which is also finite.
We consider this quantity $F_{i}$ to represent the probability of a candidate solution $i$ occurring, given the space of all possible solutions.
We further show in the discrete and continuous space that for large $N$, the candidate distribution maximising $S$ occurs with probability 1.

Consider the analogy of having a physical system made up of particles.
The different ways, $W$, in which we can organise this system of $N$ particles with $T$ distinguishable groups each containing $n_{t}$ particles, can be expressed as the combinatorial

\begin{equation}
 W = \frac{N!}{\prod_{t=1}^{T}n_{t}!} ,
\end{equation}
where $\sum_{t}n_{t}=N$. If we consider the logarithm of the above term, we can invoke Stirling's approximation that
$\log(n!) \approx n\log(n)-n$, which is exact in the limits $N\rightarrow\infty$ and $n_{i} \rightarrow \infty$.
Using this relation, we obtain
\begin{equation}
 \log W = N\bigg(-\sum_{t=1}^{T}\frac{n_{t}}{N}\log\bigg[\frac{n_{t}}{N}\bigg]\bigg) = NS,
\end{equation}
where $S$ is the Boltzmann-Shannon-Gibbs entropy and in the continous case we identify $p(t) = n_{t}/N$, where $p(t)$ represents the probability of being in group $t$. Hence, $W = \exp(NS)$.

The number of formulations, $W_{\text{maxent}}$, in which the maximum entropy realisation is more probable than any other realisation can be succinctly expressed as
\begin{equation}
 \frac{W_{\text{maxent}}}{W_{\text{other}}} =  \exp\bigl(N(S_{\text{maxent}}-S_{\text{other}})\bigr)
\end{equation}
which exposes that in the limit of large $N$, the maximum entropy solution dominates other solutions satisfying the same constraints. 
More significantly, we can also show that it dominates the space of \textit{all} solutions.
Let $\sum_{i}W_{i}$ denote the total number of ways in which the system can be configured for all possible underlying densities satisfying the constraints. 

If we consider the ratio between this term and the number of ways the maximum entropy distribution can be configured, we observe that
\begin{equation}
\label{maxentdomination}
\frac{W_{\text{maxent}}}{\sum_{i}W_{i}} = \frac{\exp(NS_{\text{max}})}{\sum_{i}\exp(NS_{i})} = \frac{1}{\sum_{i}\exp \left(N(S_{i}-S_{\text{max}})\right)} \xrightarrow{N \rightarrow \infty}{} 1
\end{equation}
where we have exploited the fact that one of the $S_{i}$ is $S_{\text{max}}$ and that $S_{i\neq j} < S_\text{max}$. 
More formally, we consider the probability mass about maxima in the functional describing all possible configurations, which is characterised via their entropy, $S$:
\begin{equation}
 W_{\text{total}} = \int \exp(NS)[\mathcal{D}P] = \idotsint_{-\infty}^{\infty} \exp(NS) \Pi_{x}dP(x)
\end{equation}
When $N\rightarrow\infty$, the maximum value of $S$ accounts for the majority of the integral's mass.  

To see this consider the ratio of functional integrals,
	\begin{equation}
	 \frac{W_{\text{maxent}}}{W_{\text{total}}} = \frac{\int \exp(NS_{\text{max}}) dP_{\text{maxent}}(x)}{\idotsint_{-\infty}^{\infty} \exp(N(S_{\text{max}} + (S-S_{\text{max}})^{2}S''|_{S=S_{\text{max}}}/2!..))\Pi_{x}dP(x)} 
	\end{equation}
which tends to 1 as $N \rightarrow \infty$. The argument is the continuous version of that displayed in Equation \eqref{maxentdomination}, where we have used Laplace's method for the multivariate case, and the definition of a probability density and the functional integral.

Laplace's method gives a theoretical basis for the canonical distributions in statistical mechanics. Its equivalent in the complex space, the method of steepest descent, in Feynman's path integral formulation, shows that points in the vicinity of the extrema of the action functional (the classical mechanical solution), contribute maximally to the path integral. This is the corresponding result for the matrix eigenvalue distributions. 

\subsection{Validation on Synthetic Data}

\begin{figure}[t]
\label{medianmatrices}
\centering 
\includegraphics[width=0.75\textwidth]{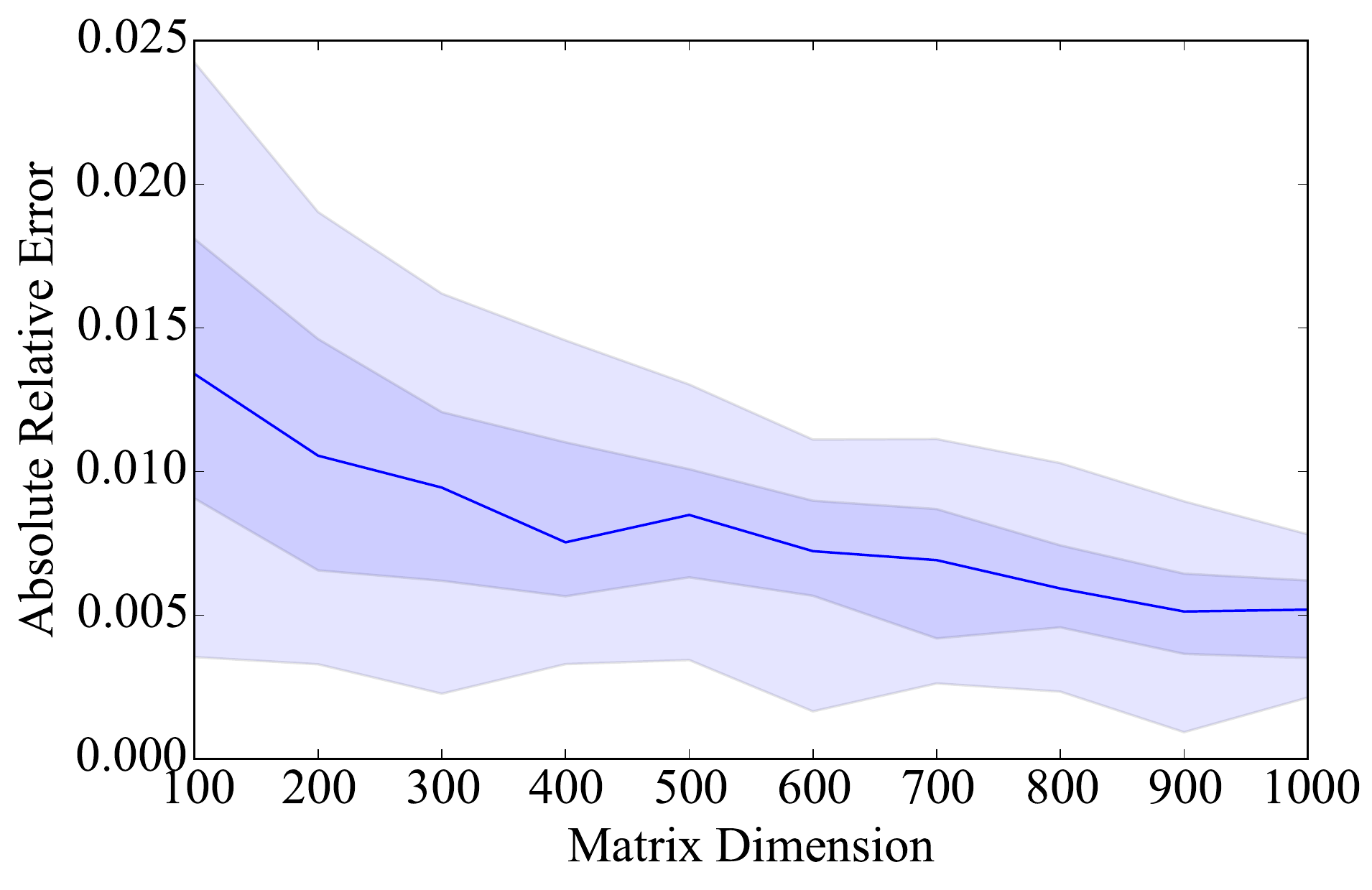}
\caption{The absolute relative error with respect to matrix dimensionality. Plotted is the median error, with 30-70 and 10-90 percentile region shaded in dark and light blue respectively.}\label{fig:error_n}
\end{figure}

We generate random, diagonally dominant positive semi-definitive matrices, $M$, which are constructed as
 \begin{equation}
 M = \frac{A^\top A}{||A^\top A||_2} + I_{N},
 \end{equation}
where $A \in \mathcal{R}^{N \times N}$ is an $N \times N$ matrix filled with Gaussian random variables and $I$ is the identity.
In order to test the hypothesis that the maximum entropy solution dominates the space of possible solutions with increasing matrix size, we investigate the relative error of the log determinant $L$, for $100 \leq N \leq 1000$. 
As can be seen in Figure \ref{fig:error_n}, there is a clear decrease in relative error for all plotted percentiles with increasing matrix size $N$.


\section{Experiments}
\label{sec:exp}

So far, we have supplemented the theoretic foundations of our proposal by devising experiments on synthetically constructed matrices.
In this section, we extend our evaluation to include real matrices obtained from a variety of problem domains, and demonstrate how the results obtained using our approach consistently outperform competing state-of-the-art approximations.
Moreover, in order to demonstrate the applicability of our method within a practical domain, we highlight the benefits of replacing the exact computation of the log determinant term appearing in the log likelihood of a Gaussian Markov random field with our maximum entropy approximation.

\subsection{UFL Datasets}

\begin{figure}[t]
 \centering 
 \includegraphics[width=\textwidth]{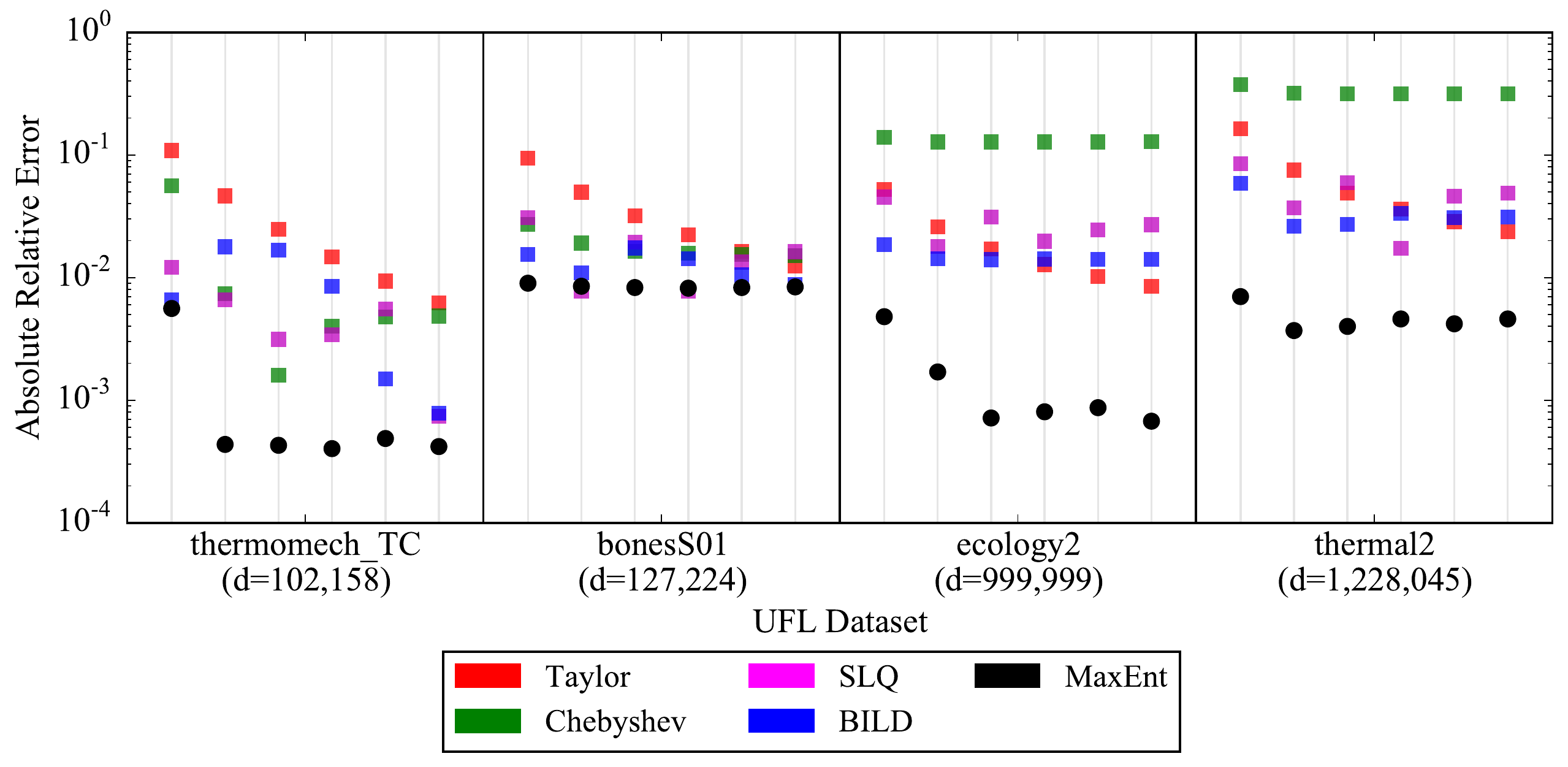}
 \caption{Comparison of competing approaches over four UFL datasets.
  Results are also shown increasing computational budgets, i.e. 5, 10, 15, 20, 25 and 30  moments respectively.
  Our method obtains substantially lower error rates across 3 out of 4 datasets, and still performs very well on 'bonesS01'.}
  \label{fig:ufl_1}
\end{figure}

While the ultimate goal of this work is to accelerate inference of large-scale machine learning algorithms burdened by the computation of the log determinant, this is a general approach which can be applied to a wide variety of application domains. 
The SuiteSparse Matrix Collection~\cite{Davis11} (commonly referred to as the set of UFL datasets) is a collection of sparse matrices obtained from various real problem domains. 
In this section, we shall consider a selection of these matrices as `matrices in the wild' for comparing our proposed algorithm against established approaches.
In this experiment we compare against Taylor~\cite{Aune2014} and Chebyshev~\cite{Han2015} approximations, stochastic lanczos quadrature (SLQ)~\cite{Ubaru2016} and Bayesian inference of log determinants (BILD)~\cite{Fitzsimons2017}.
In Figure~\ref{fig:ufl_1},  we report the absolute relative error of the approximated log determinant for each of the competing approaches over four different UFL datasets.
Following~\cite{Fitzsimons2017}, we assess the performance of each method for an increasing computational budget, in terms of matrix vector multiplications, which in this case corresponds to the number of moments considered.
It can be immediately observed that our entropic approach vastly outperforms the competing techniques across all datasets, and for any given computational budget.
The overall accuracy also appears to consistently improve when more moments are considered.

\begin{table}
\caption{
	Comparison of competing approximations to the log determinant over additional sparse UFL datasets.
	The technique yielding the lowest relative error is highlighted in bold, and our approach is consistently superior to the alternatives.
	Approximations are computed using 10 moments estimated with 30 probing vectors.
}

\begin{center}
	\begin{tabular}{lcccccc}
	\toprule
	\emph{Dataset} & Dimension &Taylor & Chebyshev & SLQ & BILD & MaxEnt\\
	\midrule
	shallow\_water1 & 81,920 & \textbf{0.0023}&  0.7255&  0.0058& 0.0163 & 0.0030\\
	shallow\_water2 & 81,920 & 0.5853&  0.9846&  0.9385&  1.1054& \textbf{0.0051}\\
	apache1 & 80,800 & 0.4335&  0.0196&  0.4200&  0.1117& \textbf{0.0057}\\
	finan512 & 74,752 & 0.1806 & 0.1158 & 0.0142 &  \textbf{0.0005}& 0.0171\\
	obstclae & 40,000 &0.0503 & 0.5269 & 0.0423 & 0.0733 & \textbf{0.0026}\\
	jnlbrng1 & 40,000 &0.1084 & 0.2079 & 0.0465 & 0.0805 & \textbf{0.0158}\\
	\bottomrule
	\end{tabular}
	\end{center}
\label{tab:ufl_larger}
\end{table}

Complementing the previous experiment, Table~\ref{tab:ufl_larger} provides a further comparison on a range of other sample matrices which are large, yet whose determinants can be computed by standard machines in reasonable time (by virtue of being sparse). 
For this experiment, we consider 10 estimated moments using 30 probing vectors, and their results are reported for the aforementioned techniques. 
The results presented in Table \ref{tab:ufl_larger} are the relative error of the log determinants \textit{after} they have been normalised using Gershgorin intervals~\cite{Gershgorin1931}.
We note, however, that the methods improve at different rates as more raw moments are taken.

\subsection{Computation of GMRF Likelihoods}

Gaussian Markov random fields (GMRFs)~\cite{Rue2005} specify spatial dependence between nodes of a graph with Markov properties, where each node denotes a random variable belonging to a multivariate joint Gaussian distribution defined over the graph.
These models appear in a wide variety of applications, ranging from interpolation of spatio-temporal data to computer vision and information retrieval. 
While we refer the reader to~\cite{Rue2005} for a more comprehensive review of GMRFs, we highlight the fact that the model relies on a positive-definite precision matrix $Q_\theta$ parameterised by $\theta$, which defines the relationship between connected nodes; given that not all nodes in the graph are connected, we can generally expect this matrix to be sparse.
Nonetheless, parameter optimisation of a GMRF requires maximising the following equation:
$$
\log p(\xvect \mid \theta) = \frac{1}{2}\log\bigl(\Det(Q_\theta)\bigr) - \frac{1}{2}\xvect^\top Q_\theta\xvect - \frac{n}{2}\log(2\pi) 
$$
where computing the log determinant poses a computational bottleneck, even where $Q_{\theta}$ is sparse.
This arises because it is possible for the Cholesky decomposition of a sparse matrix with zeros outside a band of size $k$ to be nonetheless dense \textit{within} that bound.
Thus, the Cholesky decomposition is still expensive to compute.

\begin{figure}[t!]
	\centering
	\includegraphics[width=.8\textwidth]{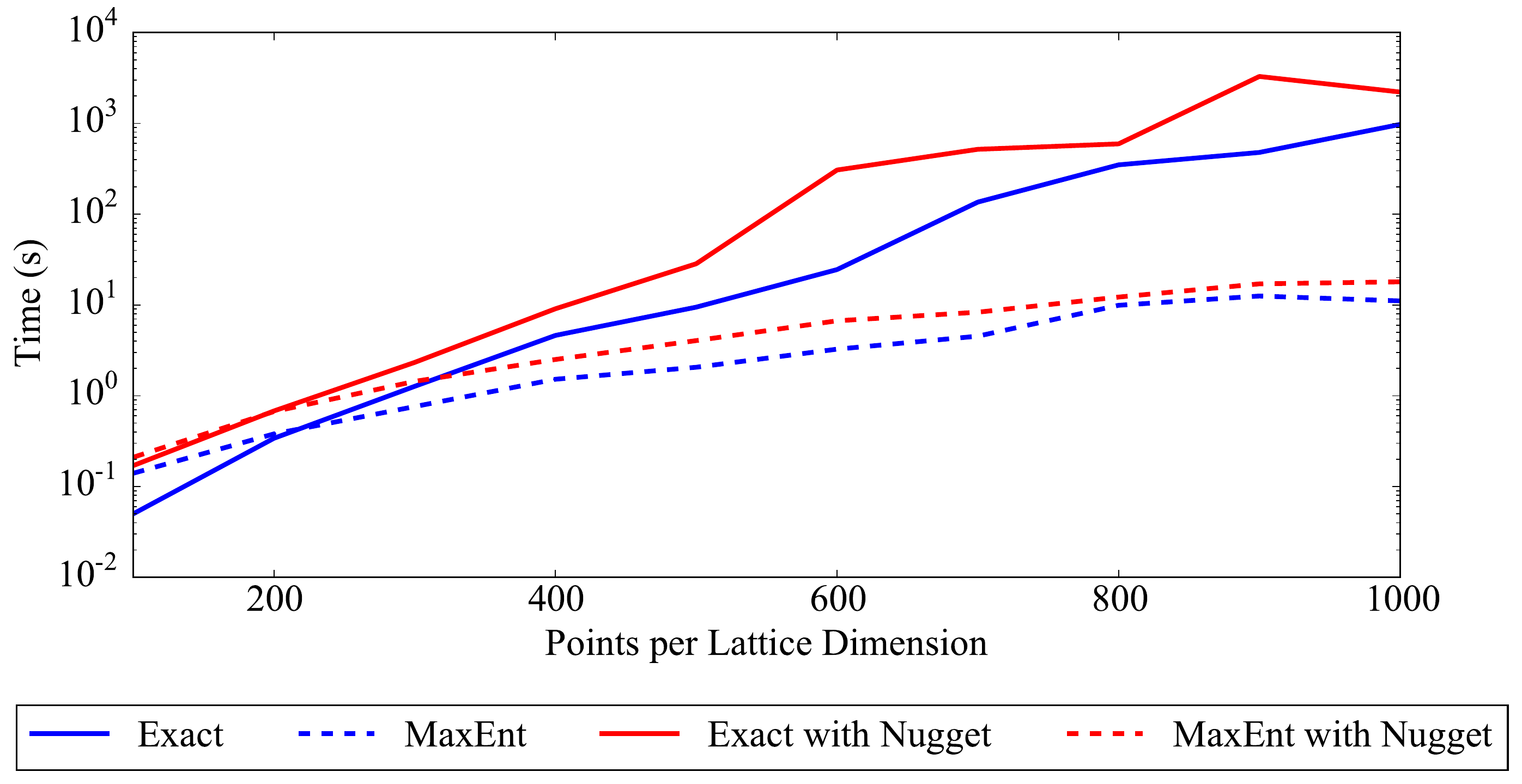}
	\caption{Time in seconds for computing the log likelihood of a GMRF via Cholesky decomposition or using our proposed MaxEnt approach for estimating the log determinant term.
	Results are shown for GMRFs constructed on square lattices with increasing dimensionality, with and without a nugget term.
	}
	\label{fig:gmrf}
\end{figure}

Following the experimental set-up and code provided in~\cite{Guinness2015}, in this experiment we evaluate how incorporating our approximation into the log likelihood term of a GMRF improves scalability when dealing with large matrices, while still maintaining precision.
In particular, we construct lattices of increasing dimensionality and in each case measure the time taken to compute the log likelihood term using both approaches.
The precision kernel is parameterised by $\kappa$ and $\tau$ \cite{lindgren2011explicit}, and is explicitly linked to the spectral density of the Mat\'ern covariance function for a given smoothness parameter.
We repeat this evaluation for the case where a nugget term, which denotes the variance of the non-spatial error, in included in the constructed GMRF model.
Note that for the maximum entropy approach we employ 30 sample vectors in the stochastic trace estimation procedure, and consider 10 moments.
As illustrated in Figure~\ref{fig:gmrf}, the computation of the log likelihood is orders of magnitude faster when computing the log determinant using our proposed maximum entropy approach.
In line with our expectations, this speed-up is particularly significant for larger matrices.
Similar improvements are observed when a nugget term is included.
Note that we set $\kappa = 0.1$ and $\tau = 1$ for this experiment.

\begin{figure}[t!]
	\centering
	\includegraphics[width=.8\textwidth]{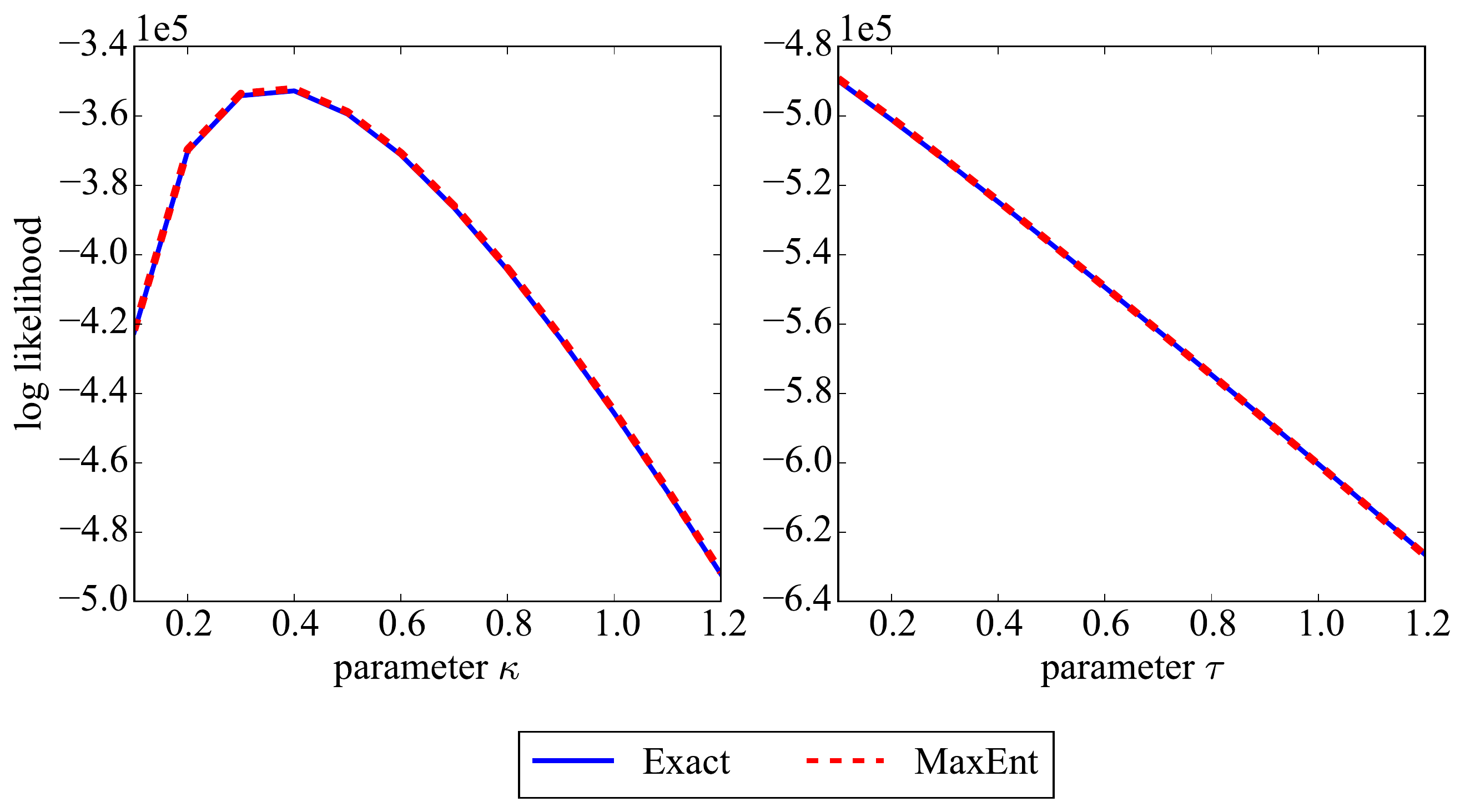}
	\caption{The above plots indicate the difference of log likelihood between exact computation of the likelihood and the maximum entropy approach for a range of hyperparameters of the model.
	We note that the extrema of both exact and approximate inference align and it is difficult to distinguish the two lines.
	}
	\label{fig:gmrf_param_sweep}
\end{figure}

Needless to say, improvements in computation time mean little if the quality of inference degrades. 
Figure \ref{fig:gmrf_param_sweep} illustrates the comparable quality of the log likelihood for a various settings of $\kappa$ and $\tau$, and the results confirm that our method enables faster inference without compromising on performance. 

\section{Conclusion}
Inspired by the probabilistic interpretation introduced in~\cite{Fitzsimons2017}, in this work we have developed a novel approximation to the log determinant which is rooted in information theory.
While lacking the uncertainty quantification inherent to the aforementioned technique,  this formulation is appealing because it uses a comparatively less informative prior on the distribution of eigenvalues, and we have also demonstrated that the method is theoretically expected to yield superior approximations for matrices of very large dimensionality.
This is especially significant given that the primary scope for undertaking this work was to accelerate the log determinant computation in large-scale inference problems.
As illustrated in the experimental section, the proposed approach consistently outperforms all other state-of-the-art approximations by a sizeable margin.

Future work will include incorporating the empirical Monte Carlo variance of the stochastic trace estimates into the inference scheme, extending the method of maximum entropy to include noisy constraints, and explicitly evaluating the ratio of the functional integrals for large matrices to obtain uncertainty estimates similar to those in \cite{Fitzsimons2017}. 
We hope that the combination of these advancements will allow for an apt active sampling procedure given pre-specified computational budgets.

\subsubsection{Acknowledgements}
Part of this work was supported by the Royal Academy of Engineering and the Oxford-Man Institute. MF gratefully acknowledges support from the AXA Research Fund. 


\bibliographystyle{splncs03}
\bibliography{references}

\end{document}